\documentclass[12pt]{article}

\usepackage{graphicx}
\usepackage{epsfig}
\usepackage{epstopdf}

\usepackage{amssymb, amsmath, amsthm}
\usepackage{float}

\newtheorem{dfn}{Definition}
\newtheorem{defn}[dfn]{Definition}

\newtheorem{rem}[dfn]{Remark}

\newtheorem{thm}{Theorem}
\newtheorem{lem}{Lemma}

\newtheorem{cor}[thm]{Corollary}

\newtheorem*{spec}{Special case}

\def\proof{\par\medskip\noindent{\it Proof: }}

\def\>{\rangle}
\def\<{\langle}

\def\3{\ss}
\def\8{\infty}

\begin{document}

\title{The surface complex of Seifert fibered spaces}

\author{Jennifer Schultens}

\maketitle

\begin{abstract}
We define the surface complex for $3$-manifolds and embark on a case study in the arena of Seifert fibered spaces.  The base orbifold of a Seifert fibered space captures some of the topology of the Seifert fibered space, so, not surprisingly, the surface complex of a Seifert fibered space always contains a subcomplex isomorphic to the curve complex of the base orbifold.  
\end{abstract}
 
\section{Introduction}

The curve complex associated with a surface continues to generate interesting research, both as an object of study in its own right and as a tool for understanding the mapping class groups of surfaces and Heegaard splittings of $3$-manifolds.  It is defined in terms of disjoint submanifolds of a surface.  Building on the success of the curve complex, we move the discussion up one dimension in an attempt to identify and describe the surface complex associated with a $3$-manifold.  Though the complex studied in \cite{Schultens2018} is quite different from the surface complex, some of the reasoning applies to the surface complex of a Seifert fibered space.  Given a Seifert fibered space $M$ with base orbifold $Q$ and natural projection $p: M \rightarrow Q$, we denote the surface obtained from the base orbifold by deleting regular neighborhoods of the cone points by $\hat Q$.  Our surfaces and Seifert fibered spaces will always be connected.  We prove the following theorems:

\begin{thm} \label{nonzero}
If $M$ is a totally orientable Seifert fibered space with nonzero Euler number, then ${\cal S}(M)$ is isomorphic to the curve complex of $\hat Q$.  \end{thm}

\begin{thm} \label{sphericalbase}
If $M$ is a totally orientable Seifert fibered space with base orbifold of genus $0$ and Euler number $0$, then ${\cal S}(M)$ contains a subcomplex isomorphic to the curve complex of the surface obtained from $\hat Q$.  Moreover, ${\cal S}(M)$ is contained in the cone on this subcomplex.
\end{thm}

\begin{spec}
If $M$ is a totally orientable Seifert fibered space with base orbifold of genus $0$, Euler number $0$, and either $4$ or $5$ exceptional fibers with identical invariants, then ${\cal S}(M)$ is isomorphic to the cone on the curve complex of $\hat Q$.  \end{spec}

\begin{thm} \label{atmostd}
If $M$ is a totally orientable Seifert fibered space with Euler number $0$ and base orbifold of positive genus, then ${\cal S}(M)$ contains a subcomplex isomorphic to the curve complex of the surface obtained from $\hat Q$.  Moreover, ${\cal S}_d(M)$ is connected, for $d$ the least common multiple of $\alpha_1, \dots, \alpha_k$.  In particular, ${\cal S}(M) = {\cal S}_d(M)$.  
\end{thm}

\section{Preliminaries} \label{prelim}

For more background and notational conventions, see \cite{HatcherAT}, \cite{HempelBook}, \cite{JacoBook}, \cite{SchultensBook}, \cite{Scott}, and \cite{Thurstonnorm}.  

\begin{defn} (Surface complex)
Let $M$ be a compact orientable $3$-manifold.  We define a sequence of complexes $\{ {\cal S}_i (M)\}$, and the {\em surface complex}, ${\cal S}(M)$, of $M$ as follows:

\begin{itemize}
\item Vertices in $\{ {\cal S}_i(M)\}$ and ${\cal S}(M)$correspond to isotopy classes of compact orientable essential (incompressible, boundary incompressible and not boundary parallel) surfaces of minimal Thurston norm (properly embedded) in $M$.  
\item A pair of distinct vertices $(v_1, v_2)$ spans an edge in ${\cal S}_0(M)$ if and only if  $v_1$ and $v_2$ admit disjoint representatives.   Inductively, we construct a sequence of complexes, $\{{\cal S}_{i}(M)\}$,  whose vertices coincide, for all $i$, with those of ${\cal S}_{0}(M)$.  Given ${\cal S}_{i}(M)$, the pair of vertices $\{v_1, v_2\}$ spans an edge in ${\cal S}_{i+1}(M)$ if and only if $v_1$ and $v_2$ lie in distinct components of ${\cal S}_{i}(M)$ and admit representatives whose intersection has $i+1$ components.  

\item For all $i$, ${\cal S}_i(M)$ is a flag complex. 
\end{itemize}

The {\em surface complex} of $M$, denoted ${\cal S}(M)$, is defined to be ${\cal S}_{i_0}(M),$ where $i_0$ is the smallest natural number such that ${\cal S}_{i_0}(M)$ is connected.
\end{defn}

There are many sources on Seifert fibered spaces.  Of particular interest is H. Seifert's original paper on the subject, \cite{SeifertBook}), or W. Heil's translation, \cite{HeilTrans}.  See also \cite{Scott} and \cite{J-S1}.  

\begin{defn} \label{FST}
A {\em fibered solid torus} is a solid torus obtained as follows:  Given a solid cylinder  ${\mathbb D}^2 \times [0,1]$, glue ${\mathbb D} \times \{0\}$ 
to ${\mathbb D} \times \{1\}$ after a rotation by a rational multiple of $2\pi.$  More specifically, the rotation is by $ \frac{2\pi \nu}{\mu}$, where $\mu, \nu \in \mathbb{Z}$ and $g.c.d.(\mu,\nu) = 1$.  The resulting solid torus is denoted by $V(\nu,\mu).$  By convention, we require that $0 < \nu < \mu.$  
A {\em fiber} is a simple closed curve resulting by identification of the endpoints of intervals of the form $\{y\} \times [0, 1], y \in {\mathbb D}^2$.

An {\em exceptional} fiber is a fiber resulting by identifying $\{(0,0)\} \in {\mathbb D}^2 \times [0,1]$ to $\{(0,1)\} \in {\mathbb D}^2 \times [0,1]$.  All other fibers are {\em regular} fibers. 
\end{defn}

\begin{figure}[h]
\includegraphics[scale=.4]{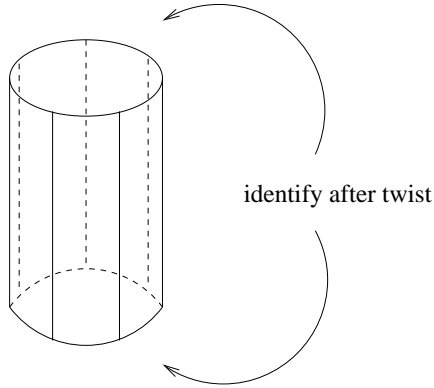}
\centering
\caption{\sl A fibered solid torus}
\label{fibtor}
\end{figure}

\begin{defn}\label{DSFS}
A {\em Seifert fibered space} $M$ is a compact connected $3$-manifold that admits a decomposition into
disjoint circles each of which has a neighborhood that is 
homeomorphic to a fibered solid torus via a homeomorphism that takes circles to fibers.  The circles into which $M$ falls are called {\em fibers} of $M.$ A particular decomposition of $M$ into fibers is called a {\em fibration} of $M.$  A fiber of $M$ is {\em exceptional} if it is exceptional fiber in a fibered solid torus neighborhood and {\em regular} otherwise.
\end{defn}

\begin{defn} \label{orbit}
Given a Seifert fibered space $M$ and a fibration, we form a quotient space, $Q$, by identifying each fiber to a point.   The quotient map is denoted by $p: M \rightarrow Q.$  Topologically, the quotient space is a surface.  However, if nearby regular fibers wrap around the exceptional fiber $e$ exactly $\mu$ times, then we declare $p(e)$ to be a {\em cone point} of $Q$ of {\em multiplicity} $\mu.$ Thus $Q$  is an orbifold, called the {\em base orbifold}.   

If the underlying surface of $Q$ is a sphere, then we say that $M$ has {\em spherical} base orbifold.  
We denote the surface obtained from $Q$ by removing regular neighborhoods of the cone points by $\hat Q$.  
\end{defn}

\begin{figure}[h]
\includegraphics[scale=.6]{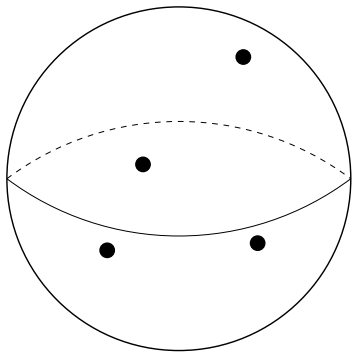}
\centering
\caption{\sl A spherical base orbifold with $4$ cone points}
\label{BaseOrbifold}
\end{figure}

Note that even for an orientable Seifert fibered space,  the base orbifold can be orientable or nonorientable.  Indeed, the twisted circle bundle over the M\"{o}bius band is homeomorphic to a twisted $I$-bundle over the Klein bottle, an orientable $3$-manifold and a Seifert fibered space.  Its double is of interest, because it admits two distinct, though homeomorphic, Seifert fibrations.  Both are twisted circle bundles over the Klein bottle.

\begin{defn}
A Seifert fibered space $p: M \rightarrow Q$ is {\em totally} orientable if $M$ and $\hat Q$ are orientable.
\end{defn}

The fibering of a Seifert fibered space enables local isotopies positioning a surface to either coincide with or be transverse to the fibers.  That this is possible globally is the content of the following theorem:

\begin{thm}\cite[VI.34]{JacoBook}   \label{Jaco}
(Jaco) Let $F$ be a connected, two-sided, essential surface in an orientable Seifert fibered space $M.$ Then one of the following holds: 

(i) $F$ is non separating in $M$ and is a fiber in a fibration of $M$ as a surface bundle over the circle;

(ii) $F$ separates $M$ and $M = M_1 \cup M_2,$ where $\partial M_i = F$ and $M_i$ is a twisted $I$-bundle over a compact surface;

(iii) $F$ is an annulus or torus and $F$ is saturated, {\em i.e.,} consists of fibers, in some Seifert fibration of $M.$
 
 \end{thm}

In Case (i), $M$ can be described both as a Seifert fibered space and as a surface bundle over the circle.  Moreover, $M$ is obtained by identifying the boundary components of $F \times I$.  The quotient map $p: M \rightarrow Q$ restricts to $F$ to give an irregular covering $p|_F \rightarrow Q$.  The fibers of $M$ consist of unions of intervals of the form $\{point\} \times I$.  In the case where $F$ is a torus, and in a handful of other cases, the Seifert fibered space will admit more than one Seifert fibration.  Case (ii) does not occur if the base orbifold of $M$ is orientable.  Case (iii) describes vertical surfaces. 

\begin{rem} \label{verticalregular}
A vertical surface $S$ in a Seifert fibered space $M$ projects to a simple curve.   If $S$ contained an exceptional fiber, $e$, then $e$ would project to a point, $p$, in $Q$ on a simple curve $s$.  Nearby regular fibers on $s$ would wrap around $e$.   This makes $e$ a branch locus of order $\mu$.  For $\mu \geq 3$, this is impressible.  Hence $\mu \leq 2$.  Furthermore, if $\mu = 2$, then $S$ meets a fibered solid torus neighborhood of $e$ is a M\"{o}bius band, implying that $S$ is not orientable.  Thus vertical surfaces consist of regular fibers.  
\end{rem}

\begin{lem} \label{horidistinct}
If $F, F'$ are disjoint horizontal surfaces in the orientable Seifert fibered space $M$ with orientable base orbifold, then $F$ and $F'$ are isotopic.   
\end{lem}

\proof
It follows from Theorem \ref{Jaco} that $F$ and $F'$ are of type (ii) or (iii).  Case (iii) does not occur because the base orbifold is orientable.  Thus $F'$ is properly embedded in $F \times I$ and must hence be isotopic to $F$.  
\qed
 
\begin{lem} \label{horinotvert}
Let $F$ be a two-sided essential surface in an orientable Seifert fibered space  $M.$ 
If $F$ is connected and horizontal with respect to a given Seifert fibration, then it is neither isotopic nor homologous to a vertical surface with respect to this fibration.  Likewise, if $F$ is connected and vertical with respect to a given Seifert fibration, then it is neither isotopic nor homologous to a horizontal surface with respect to this fibration.
\end{lem}
 
This follows from the proofs of Theorems VI.26 and VI.34 in \cite{JacoBook}.  

Horizontal surfaces behave quite differently than vertical surfaces, but here too, the considerations in \cite{Schultens2018} provide key insights.
Recall the following theorem, see for instance \cite[Corollary 5]{Schultens2018}, about the correspondence between horizontal surfaces and $2$-dimensional homology classes:

\begin{thm} \label{hori}
If $F$ and $F'$ are homologous horizontal surfaces in an orientable Seifert fibered space, then $F$ and $F'$ are isotopic. 
\end{thm}

By definition, every Seifert fibered space is foliated by circles.  One must not confuse this Seifert fibration with an honest fibering as a circle bundle over the surface or as a surface bundle over the circle.  Of course, some Seifert fibered spaces, such as the manifolds $(surface) \times {\mathbb S}^1$, exhibit both a (trivial) Seifert fibering, a fibering as a circle bundle over a surface, and a fibering as a surface bundle over the circle.
One of the invariants of a Seifert fibered space $M$ is the Euler number.   It measures the obstruction of $M$ to being a surface bundle over the circle.  Only those Seifert fibered spaces with Euler number $0$ fiber as surface bundles over the circle. 

More specifically, 
let $M$ be Seifert fibered space with base orbifold $Q$ and projection $p: M \rightarrow Q$.  Consider the Seifert fibered space $M*$ with base orbifold $Q*$ obtained from $M$ by removing a fibered solid torus neighborhood $T$ of a regular fiber.  Then $M*$ admits a horizontal surface $\tilde Q*$.  The projection map restricts to an orbifold covering of $p|_{\tilde Q*}: \tilde Q* \rightarrow Q*$.  In general, $\tilde Q*$ will not extend to a horizontal surface $\tilde Q$ covering $Q$, the obstruction being an integer $b$, the oriented intersection number of a meridian of $T$ with $\partial Q$.  In other words, the horizontal surface $\tilde Q*$ extends to a horizontal surface $\tilde Q$ in $M$ if and only if $b = 0$.  

Seifert fibered spaces provide a superb arena in which to develop tools for studying $3$-manifolds.  The existence of a quotient map allows much of the information needed to identify a Seifert fibered space to be captured in the base orbifold.  This is true of algebraic information, for instance in our computations of the fundamental group.  To a certain degree, it allows us to understand incompressible surfaces in Seifert fibered spaces.
A Seifert fibered space $M$ is completely determined by a set of invariants computed from the base orbifold, an invariant called the {\em Euler number}, and the $\mu$s and $\nu$s of the exceptional fibers.  In the case of a closed totally orientable Seifert fibered space the invariants are as follows:

$$M = < g, b, (\alpha_1, \beta_1), \dots, (\alpha_k, \beta_k) >$$

The number $g$ is the genus of the base orbifold $Q$ of $M$.  The number $b$ is the Euler number.  The pair $(\alpha_i, \beta_i)$ provides the slope of the meridian disk of a solid torus neighborhood of the $i$th exceptional fiber in terms of ``external" coordinates.  As in Definition \ref{FST}, the pair $(\mu_i, \nu_i)$ provides ``internal" coordinates describing the fibering of  fibered solid torus neighborhood $V_i$ of the $i$th exceptional fiber.  On the other hand, denote the surface obtained by deleting a regular neighborhood of a point from $\hat Q$ by $\hat Q^-$.  Then $\hat Q^-$ is a surface with boundary components $(c_1, \dots, c_{k+1})$.  For each $i$, choose a point $p_i$ in $c_i$.  Then the pair $(c_i \times \{p_i\}, \{p_i\} \times {\mathbb S}^1)$ of curves in the boundary component $c_i \times {\mathbb S}^1$ of $\hat Q^- \times {\mathbb S}^1$ provides coordinates in which to express the slope of the meridian disk of the $i$th exceptional fiber.  We take $\alpha_i$ to be the number of times a meridian of $V_i$ meets $\{p_i\} \times {\mathbb S}^1$ (or any regular fiber, hence $\alpha_i = \mu_i$ ) and $\beta_i$ to be the number of times it meets $c_i \times \{p_i\}$.  To construct $M$ from $\hat Q^- \times {\mathbb S}^1$, perform Dehn fillings of slope $(\alpha_i, \beta_i)$ along $c_i \times {\mathbb S}^1$ for $i = 1, \dots, k$, then perform a Dehn filling of slope $(1, b)$ along the remaining boundary component.

The fundamental group of $M$ can be computed from this set of invariants.  

$$\pi_1(M) = < a_1, b_1, \dots, a_g, b_g, x_1, \dots, x_k, h \; | $$
$$\; h^{-b}\Pi_1^g [a_i, b_i] \Pi_1^k x_i, [a_1, h], [b_1, h], \dots, [a_g, h], [b_g, h], [x_1, h], \dots, [x_k, h], x_1^{\alpha_1}h^{\beta_1}, \dots, x_k^{\alpha_k}h^{\beta_k} > $$

Since a Seifert fibered space admits horizontal surfaces if and only if its Euler number is $0$ and since horizontal surfaces are in 1-1 correspondence with $2$-dimensional homology classes, we compute, for $M$ a Seifert fibered space with Euler number $0$:

$$\pi_1(M) = < a_1, b_1, \dots, a_g, b_g, x_1, \dots, x_n, h \; | $$
$$\; \Pi_1^g [a_i, b_i] \Pi_1^n x_i, [a_1, h], [b_1, h], \dots, [a_g, h], [b_g, h], [x_1, h], \dots, [x_n, h], x_1^{\alpha_1}h^{\beta_1}, \dots, x_n^{\alpha_n}h^{\beta_n} > $$
Abelianizing, we obtain:
$$H_1(M) = < a_1, b_1, \dots, a_g, b_g, x_1, \dots, x_n, h \; | \;  x_1 + \dots + x_n, \alpha_1 x_1 + \beta_1 h, \dots, \alpha_n x_n + \beta_n h > $$
Relations of the form $\alpha_i x_i + \beta_i h$ yield relations between the $x_i$s.  Specifically, we obtain, for all $i, j$, the following:
$$\beta_i \alpha_j x_j + \beta_i \beta_j h = \beta_j \alpha_i x_i + \beta_j \beta_i h$$
Consider, for instance, the case $\alpha_1 = 3, \beta_1 = 2, \alpha_2 = 5, \beta_2 = 3$.  
$$9x_1 + 6h = 10x_2 + 6h$$
$$9(x_1 - x_2) = x_2$$
This allows us to compute the abelian group expllcitly:  $$< x_1, x_2 \; | \; 9x_1 = 10x_2\;> = < x_1 - x_2, x_2 \; | \; 9(x_1 - x_2) = x_2 > $$ $$= <x_1 - x_2> = {\mathbb Z}$$ The substitution used is an example of a standard procedure involving Nielsen equivalence and the Euclidean algorithm.  In general, Nielsen equivalence provides a method for reducing the number of generators.  This allows us to compute $H_1$ explicitly:
$$H_1(M) = < a_1, b_1, \dots, a_g, b_g, \eta > $$

If $M$ is a product, then $\eta = h$, the homology class of a regular fiber.  If there are exceptional fibers, $\eta$ will be represented by a simple closed curve, certainly not a fiber, that can be thought of as a Poincar\'{e} dual to a horizontal surface.  

By the Universal Coefficient Theorem and Poincar\'{e} duality, we obtain a presentation of the $2$-dimensional homology as a free group on $2g + 1$ generators.  The generators can be thought of explicitly:  The generator of $H_2(M)$ corresponding to $a_1$ will be $[\beta_1 \times {\mathbb S}^1]$ (for $\beta_1$ an embedded curve such that $[\beta_1] = b_1$), the generator of $H_2(M)$ corresponding to $b_1$ will be $[\alpha_1 \times {\mathbb S}^1]$ (for $\alpha_1$ an embedded curve such that $[\alpha_1] = a_1$), $\dots,$ the generator of $H_2(M)$ corresponding to $a_g$ will be $[b_g \times {\mathbb S}^1]$ (for $\beta_g$ an embedded curve such that $[\beta_g] = b_g$), the generator of $H_2(M)$ corresponding to $b_g$ will be $[a_g \times {\mathbb S}^1]$ (for $\alpha_g$ an embedded curve such that $[\alpha_g] = a_g$), the generator corresponding to $\eta$ will be the homology class of a horizontal surface.   Abusing notation, we will use the same notation for the generators of $H_2(M) = H^1(M) = H_1(M)$:  $a_1, b_1, \dots, a_g, b_g, \eta  $.


Intuitively speaking, exceptional fibers can ``get in the way of horizontal surfaces".  More rigorously put, a Seifert fibered space admits horizontal surfaces if and only if its Euler number is $0$.  The case of Seifert fibered spaces with spherical base turns out to be rather constrained.  In view of our knowledge of the curve complex, the theorem below tells us exactly what to expect.  

\begin{rem} \label{connectedboundary}
Analogous computations can be made for Seifert fibered spaces with boundary.  For instance, a Seifert fibered space fibered over the disk with $k$ exceptional fibers will have second relative homology equal to ${\mathbb Z}$.  
\end{rem}

\section{Finegold's torus complex} \label{torussection}

In her dissertation, Brie Finegold defines and studies both oriented and unoriented torus complexes for special linear groups over rings in dimension $n \geq 2$.  She is interested in generating systems of groups, but  the discussion is interesting in its own right.  For the group ${\mathbb Z}$, Finegold's construction of ${\cal C}({\mathbb T}^n)$ reduces to the following:

\begin{defn} (Finegold's unoriented torus complex for the integers)
\begin{itemize}  
\item Vertices of ${\cal C}({\mathbb T}^n)$ correspond to primitive, nonzero, integral vectors, $(x_1, \dots, x_n)$ $\in {\mathbb Z}^n$, up to sign.  ({\it I. e., elements of ${\mathbb QP}^{n-1}$.}) We write $[(x_1, \dots, x_n)]$ for the equivalence class of $(x_1, \dots, x_n)$.    
\item For $1 \leq k \leq n - 1$, the vertices $[v_0], [v_1], \dots, [v_k]$ span a $k$-simplex if and only if the $n \times k$ matrix $M(v_0, v_1, \dots, v_k)$ is a submatrix of an element of $SL(n, {\mathbb Z}).$  A set of $n+1$ vertices spans an $n$-simplex if and only if every subset obtained by omitting a vertex spans an $(n-1)$-simplex.  
\end{itemize}
\end{defn}

For $n= 2$, the result is the Farey complex.  The Farey complex is the curve complex of a torus.  It is important to note that the $2$-torus does not admit disjoint non isotopic essential curves.  The edges of the curve complex are therefore defined to be pairs of isotopy classes of curves that meet (transversely) in exactly one point. 

The first cohomology of ${\mathbb T}^n$ is ${\mathbb Z}^n.$ Given a $1$-dimensional cohomology class ${\bf v}$, its dual is an $(n-1)$-dimensional homology class.  Since nonzero $(n-1)$-dimensional homology classes are represented by closed oriented hypersurfaces (see Proposition 1.7.16 in \cite{MartelliBook}), a $1$-dimensional cohomology class thus corresponds to a $(n-1)$-dimensional hypersurface.   Specifically, the $1$-dimensional cohomology class ${\bf v}$ is the dual of $[S_{\bf v}]$, where $S_{\bf v}$ is the hypersurface corresponding to ${\bf v}$.  

The cohomology class  ${\bf v}$ is realized by intersection number number with $S_{\bf v}$.   It follows that, when ${\bf v}$ is primitive, the homology class $[S_{\bf v}]$ is also primitive.  It is known that such a homology class can be realized by a connected closed oriented hypersurface.  To simplify our discourse, we will say that ${\bf v}$ {\em represents} $S_{\bf v}$.  
In ${\mathbb T}^n$ a connected closed oriented hypersurface is homologous to a canonical (flat) essential oriented $(n-1)$-torus.  Hence, in the case $n = 3$, the vertices of Finegold's torus complex are in 1-1 correspondence with the vertices of ${\cal S}_0({\mathbb T}^n)$ and thus coincide with the vertices of ${\cal S}({\mathbb T}^n)$.   

Non isotopic essential tori in ${\mathbb T}^3$ necessarily meet, hence ${\cal S}_0({\mathbb T}^3)$ consists of isolated vertices.  
Recall that the mapping class group of ${\mathbb T}^3$ is $GL(3, {\mathbb Z})$.  Two vertices span an edge in ${\cal S}_1({\mathbb T}^3)$ if and only if  they are represented by tori, $T_a, T_b$, that can be isotoped to meet in a single component.  This happens if and only if there is an element of $SL( 3, {\mathbb Z})$ that takes the two vectors ${\bf a}, {\bf b}$ representing the two tori to two coordinate vectors, say $(1, 0, 0)$ and $(0, 1, 0)$.   The first two columns, ${\bf a}, {\bf b}$, of $A^{-1}$ are the vectors representing $T_a, T_b$.  Since $A^{-1} \in SL(3, {\mathbb Z})$, the two vertices ${\bf a}, {\bf b}$ span an edge in Finegold's torus complex.
Thus the $1$-skeleton of Finegold's torus complex coincides with the $1$-skeleton of ${\cal S}_1({\mathbb T}^3)$.  
         
Recall that the graph distance between vertices in a complex is the least number of edges required in an edge path between the two vertices.  This provides a metric on the set of vertices.  Among other things, Finegold establishes the following:  

\begin{thm} \label{Finegold1}
The diameter of ${\cal C}({\mathbb T}^3)$ is $2$.  
\end{thm}

The proof is elementary.  We include it here.

\proof
Let $T$ be an essential torus in ${\mathbb T}^3$.  Then the vertex $[T]$ corresponds to a primitive, nonzero, integral vector $(a, b, c)$.
The Euclidean algorithm furnishes numbers $(x, y)$ such that $ay - bx = gcd(a, b)$.  Note that, while the pairs of numbers satisfying this equation are not unique, the pair $(x,y)$ furnished by the Euclidean algorithm is unique and $gcd(x,y) = 1$.  In particular, since $(a, b, c)$ is primitive, $gcd(gcd(a,b), c) = 1$, so there are relatively prime numbers $\alpha, \beta$ such that $\alpha gcd(a, b) + \beta c = 1$.
Furthermore, since $gcd(x,y) = 1$, there are relatively prime numbers $\gamma, \delta$ such that $\gamma x + \delta y = 1$.  This shows the following: 

\[ 
\begin{vmatrix}
a & x & - \beta \delta  \\ 
b & y & \beta \gamma  \\ 
c & 0 & \alpha  \notag
\end{vmatrix}
= a y \alpha - b x \alpha + c ( x \beta \gamma + y \beta \delta ) = \alpha gcd(a, b) + \beta c = 1
\]

It follows that $(x,y,0)$ represents a torus $T'$ that intersects $T$ in a single simple closed curve, {\it i.e.,} there is an edge $((a, b, c), (x, y, 0))$.  In turn,

\[ 
\begin{vmatrix}
\delta & x & 0  \\ 
-\gamma & y & 0  \\ 
0 & 0 & 1  \notag
\end{vmatrix}
=  1
\]

 Hence there is also an edge $((x,y,0), (0,0,1))$.   Therefore $[T]$ is at most distance $2$ from the vertex $(0, 0, 1)$.   
 
 Given two distinct tori  $T_a, T_b$ in ${\mathbb T}^3$, we apply a coordinate transformation that transforms, say, the vector representing $T_b$ into $(0,0,1)$, before embarking on the above computation.  It follows that $[T_a], [T_b]$ are within distance two.  
\qed

Intuitively speaking, this proof confirms that, given any torus, $T$, intersecting the ``horizontal" coordinate torus in any number of simple closed curves, there is a ``vertical" torus, $T'$, that intersects both $T$ and the ``horizontal" coordinate torus in a single simple closed curve.

\begin{cor}
The diameter of Finegold's torus complex is $2$ as is the diameter of ${\cal S}_1({\mathbb T}^3)$.  
\end{cor}

\begin{cor}
The complex ${\cal S}_1({\mathbb T}^3)$  is connected.  Hence ${\cal S}({\mathbb T}^3)  = {\cal S}_1({\mathbb T}^3)$.  
\end{cor}

Note that the $2$-skeleton of Finegold's torus complex is a subset of the $2$-skeleton of ${\cal S}({\mathbb T}^3) $.  However, it does not coincide with the $2$-skeleton of ${\cal S}({\mathbb T}^3)$.  Rather, it is a strict subset. For instance, the triple of vertices represented by the following triple of vectors spans a $2$-simplex in ${\cal S}({\mathbb T}^3)$, but not in Finegold's torus complex:

\[ 
\begin{bmatrix}
1\\ 
0  \\ 
0 \notag
\end{bmatrix}, 
\begin{bmatrix}
0  \\ 
1  \\ 
0  \notag
\end{bmatrix},
\begin{bmatrix}
1  \\ 
1  \\ 
2  \notag
\end{bmatrix}
\]

She also proves the following, more technical result, see \cite[Section 4.1]{Finegold}:

\begin{thm} \label{Finegold2}
For $n = 2, 3$, Finegold's torus complex is simply connected.
\end{thm}

\begin{thm} \label{simplyconnected}
The surface complex of the $3$-torus is simply connected.
\end{thm}

\proof
The $1$-skeleton of ${\cal S}({\mathbb T}^3)$ coincides with the $1$-skeleton of Finegold's unoriented torus complex for $n = 3$ and the $2$-skeleton of ${\cal S}({\mathbb T}^3)$ contains the $2$-skeleton of Finegold's unoriented torus complex for $n = 3$.  Hence
${\cal S}({\mathbb T}^3)$ is simply connected.  
\qed

We observe that ${\cal S}({\mathbb T}^3)$ is locally infinite.  Indeed, consider a torus $T$ represented by $(0, 0, 1)$.   Thus, for every pair $(p, q)$ of coprime integers, the torus represented by $(r, s, 0)$, where $ps+qr = \pm 1$, meets $T$ in a single essential simple closed curve.    
 
Links of edges are also infinite.  Indeed, given tori $T_z$ and $T_y$ represented by $(0,0,1)$ and $(0, 1, 0)$, together with a torus $T_r$, represented by $(1, r, 0)$, for any $r$, will represent the vertices of a $2$-simplex.  

\begin{lem}
The dimension of ${\cal S}({\mathbb T}^3) $ is at least $6$.
\end{lem}

\proof
The $7$-tuple of vertices represented by the following vectors spans a $6$-simplex:

\[ 
\begin{bmatrix}
1\\ 
0  \\ 
0 \notag
\end{bmatrix}, 
\begin{bmatrix}
0  \\ 
1  \\ 
0  \notag
\end{bmatrix},
\begin{bmatrix}
0  \\ 
0  \\ 
1  \notag
\end{bmatrix},
\begin{bmatrix}
1\\ 
1  \\ 
0 \notag
\end{bmatrix}, 
\begin{bmatrix}
1  \\ 
0  \\ 
1  \notag
\end{bmatrix},
\begin{bmatrix}
0  \\ 
1  \\ 
1  \notag
\end{bmatrix},
\begin{bmatrix}
1\\ 
1  \\ 
1 \notag
\end{bmatrix}, 
\begin{bmatrix}
1  \\ 
2  \\ 
0  \notag
\end{bmatrix}
\]
\qed

\section{Products}

Some of the geometric insights concerning ${\cal S}({\mathbb T}^3)$ extend to product manifolds.  Let $F$ be a connected closed orientable surface and $M = F \times {\mathbb S}^1$ the corresponding product manifold.  To understand the surface complex of $M$, we must understand how and when surfaces meet.  The following theorem reduces this task to considering surfaces that meet at most once:

\begin{thm} \label{atmostonce}
Let $F$ be a connected closed orientable surface.  If $M = F \times {\mathbb S}^1$, then ${\cal S}_1(M)$ is connected.  In particular, ${\cal S}(M) = {\cal S}_1(M)$.
\end{thm}

\proof
Recall that $H_2(M) = < a_1, b_1, \dots, a_g, b_g, \eta > $.  The generator $\eta$ is represented by $F \times \{point\}$.  We will show that the component of ${\cal S}_1(M)$ containing $v_{\eta} = [F \times \{point\}]$ contains all vertices of ${\cal S}_0(M)$.  A vertical surface is isotopic to $\gamma \times {\mathbb S}^1$, for $\gamma$ a simple closed curve in $F$.  Thus its intersection $$(\gamma \times {\mathbb S}^1) \cap (F \times \{point \}) = \gamma \times \{point\}$$
is connected.  Hence the vertex $[\alpha \times {\mathbb S}^1]$ spans an edge with $v_{\eta}$.  

Let $S$ be any horizontal surface.  Since $[S] = m c + n \eta$, for some primitive class $c \in <a_1, b_1, \dots, a_g, b_g,>,$ where $m, n$ are relatively prime integers, Theorem \ref{horidistinct} tells us that  $S$ is isotopic to the double curve sum of $m$ parallel copies of $\gamma \times {\mathbb S}^1$ (for an embedded curve $\gamma$ in $F$ such that $[\gamma] = c$) and $n$ parallel copies of $F \times \{point\}$.  We may assume (by deliberate choice of generators of $H_2(M)$) that $m, n \geq 0$.  

Choose an embedded curve $\gamma'$ in $F$ that meets $\gamma$ exactly once transversely.  Since $(\gamma' \times {\mathbb S}^1) \cap S$ is homologous to $m [(\gamma \cap \gamma') \times {\mathbb S}^1] + n [\gamma' \times \{point\}]$ and is an embedded multi-curve, it is, in fact, the torus link representing this homology class.  Since $m, n$ are relatively prime, it is a torus knot, in particular, it is connected.  (This can be seen explicitly by restricting the double curve sum that produces $S$ to the vertical torus $\gamma' \times {\mathbb S}^1$.) 

The vertex corresponding to $S$ thus lies in the same component of ${\cal S}_1(M)$ as the vertex corresponding to $\gamma' \times {\mathbb S}^1$.  Since $\gamma' \times {\mathbb S}^1 \cap F$ consists of a single component, the vertex corresponding to $\gamma' \times {\mathbb S}^1$ lies in the same component of ${\cal S}_1(M)$ as $v_{\eta}$.  
\qed

In fact, the proof above shows more:

\begin{cor}
If $M = F \times {\mathbb S}^1$, then ${\cal S}(M)$ has diameter at most $4$.
\end{cor}

\proof
A vertical surface meets $F \times \{point\}$ in exactly one component.  It follows that vertices represented by vertical surfaces are distance one from $v_{\eta}$.  The proof above shows that the distance between $v$ and any other horizontal surface is at most $2$.  
\qed

\section{Seifert fibered spaces}

At the outset, we consider the subcomplex ${\cal S}_0(M)$ of the surface complex of a Seifert fibered space.  As it turns out, the vertical surfaces span an interesting subcomplex of ${\cal S}_0(M)$.  The theorem below echos an analogous theorem concerning the Kakimizu complex of Seifert fibered spaces.  (See \cite{Schultens2018}.)

\begin{thm} \label{main}
Let $M$ be a totally orientable Seifert fibered space with orientable base orbifold $Q$.  Then the components of ${\cal S}_0(M)$ consist of isolated vertices represented by isotopy classes of horizontal surfaces in $M$ (if any) along with a complex that is naturally isomorphic to the curve complex of $\hat Q$.
\end{thm}

To prove Theorem \ref{main}, we need the following lemma:

\begin{lem} (\cite[Lemma 39]{Schultens2018}) \label{key}
Let $M$ be a closed totally orientable Seifert fibered space with base orbifold $Q$. Let $F$ and $F'$ be connected oriented essential vertical surfaces in $M$ and denote $p(F)$ by $b$ and $p(F')$ by $b'$.   If $F$ and $F'$ are isotopic, then $b$ and $b'$ are isotopic in $\hat Q$. \end{lem}

\proof
(Theorem \ref{main}) By Lemma \ref{horidistinct} isotopy classes of any horizontal surfaces in $M$ will correspond to isolated vertices.  By Lemma \ref{horinotvert} the isolated vertices corresponding to horizontal surfaces, if any, are distinct from vertices corresponding to vertical surfaces.  By Lemma \ref{key}, two vertical surfaces are isotopic (disjoint) if and only if their projections onto the base orbifold are isotopic (disjoint) in $\hat Q$.  Theorem \ref{main} now follows.  
\qed

\proof (Theorem \ref{nonzero})
Since there are no horizontal surfaces, Theorem \ref{nonzero} follows from Theorem \ref{main}.
\qed

\begin{thm}
If $M = F \times {\mathbb S}^1$, then ${\cal S}(M)$ contains a subcomplex isomorphic to the cone on the curve complex of $F$.
\end{thm}

\proof
The isotopy classes of vertical surfaces span a subcomplex ${\cal C}$ isomorphic to the curve complex of $F$.  The surface $F \times \{point\}$ meets each vertical surface in a single component, hence corresponds to a vertex that spans an edge with each of the vertices of ${\cal C}$. 
\qed

\proof (Theorem \ref{sphericalbase})
By Theorem \ref{main},  the components of ${\cal S}_0(M)$ consist of isolated vertices represented by isotopy classes of horizontal surfaces in $M$ (if any) along with a complex isomorphic to the curve complex of $\hat Q$.

Since $H_2(M)$ is ${\mathbb Z}$, there is, up to sign, a single primitive $2$-dimensional homology class ${\bf 1}$.  By Theorem \ref{horidistinct}, it follows that there is a unique unoriented isotopy class of horizontal surfaces representing ${\bf 1}$.   Denote this surface by $F$ and the vertex representing it by $v_{\eta}$.   Thus ${\cal S}_0(M)$ consist of the single isolated vertex $v_{\eta}$ along with a complex isomorphic to the curve complex of $\hat Q$.

The map $p$ restricts to an orbifold covering $p|_F: F \rightarrow Q$ that is unramified away from the (preimages) of cone points.  Denote the covering degree of this map by $d$.  Given a vertical torus $T$, $p$ restricts to a regular covering of degree $d$ on $F \cap T$.  Thus this intersection consists of at most $d$ components.  Since $F \cap T$ is a torus link, the specific number, $d_T$, of such components will be a divisor of $d$.  Set $$d_{min} = min\{d_T \; | \; T \; an \; essential \; vertical \; torus \; in \; M\}.$$ 

A vertex in the subcomplex of ${\cal S}_0(M)$ isomorphic to the curve complex of $\hat Q$ spans an edge with $v_{\eta}$ if and only if it is represented by a vertical torus that meets $F$ in exactly $d_{min}$ components.  Since there is at least one such vertex, by definition of $d_{min},$
${\cal S}_{d_{min}}(M)$ is connected.  In particular, ${\cal S}(M) = {\cal S}_{d_{min}}(M)$.  \qed

In the special case where there are at most $5$ exceptional fibers with identical invariants, there is exactly one essential vertical torus up to homeomorphism (though not up to isotopy).  This means that every vertical torus meets the horizontal surface in the same number, $d_{min}$, of components.  Therefore ${\cal S}_{d_{min}}(M)$, and hence also ${\cal S}(M)$, is the cone on the subcomplex isomorphic to the curve complex of $\hat Q$.  

\proof (Theorem \ref{atmostd})
By Theorem \ref{main},  the components of ${\cal S}_0(M)$ consist of isolated vertices represented by isotopy classes of horizontal surfaces in $M$ (if any) along with a complex isomorphic to the curve complex of $\hat Q$.

We need to show that ${\cal S}_d(M)$ is connected.  Denote a surface realizing $\eta$ by $F$.  The map $p$ restricts to a covering $p|_F: F \rightarrow Q$ that is unramified away from the (preimages) of cone points. Note that this restricted map must have minimal covering degree.  We denote this degree by $d$.  It is known that $d$ is the least common multiple of $\alpha_1, \dots, \alpha_k$.
  
 We will show that the component of ${\cal S}_d(M)$ containing $v_{\eta} = [F]$ contains all vertices of ${\cal S}_0(M)$.   Define $d_{min}$ as in the proof of Theorem \ref{sphericalbase}.  As in the proof of Theorem \ref{sphericalbase}, a vertex in the subcomplex of ${\cal S}_0(M)$ isomorphic to the curve complex of $\hat Q$ spans an edge with $v_{\eta}$ if and only if it is represented by a vertical torus that meets $F$ in exactly $d_{min}$ components.  Since there is at least one such vertex, by definition of $d_{min},$ the component of $v_{\eta}$ in ${\cal S}_{d_{min}}(M)$ contains this subcomplex.
   
Recall that $H_2(M) = < a_1, b_1, \dots, a_g, b_g, \eta > $.  
Let $S$ be any horizontal surface.  Since $[S] = m c + n \eta$, for some primitive class $c \in <a_1, b_1, \dots, a_g, b_g,>,$ where $m, n$ are relatively prime integers, Theorem \ref{hori} again tells us that  $S$ is isotopic to the double curve sum of $m$ parallel copies of $\gamma \times {\mathbb S}^1$ (for an embedded curve $\gamma$ in $F$ such that $[\gamma] = c$) and $n$ parallel copies of $F \times \{point\}$.  We may assume (by deliberate choice of generators of $H_2(M)$) that $m, n \geq 0$.  

In the case that $m = 0$ (or if $c$ is homologically trivial), then we must have $n = 1$ for $S$ to be connected, hence $S$ is isotopic to $F$.  
Otherwise  $m \neq 0$ and $c$ is homologically nontrivial, hence we can choose an embedded curve $\gamma'$ in $Q$ as in the proof of Theorem \ref{atmostonce} that meets $\gamma$ exactly once transversely and $\gamma' \times {\mathbb S}^1$ will be a vertical torus intersecting $S$ in $1$ component and $F$ in at most $d$ coherently oriented components. The vertex corresponding to $S$ thus lies in the same component of ${\cal S}_d(M)$ as the vertex corresponding to $\gamma' \times {\mathbb S}^1$ and as $v_{\eta}$.    Thus ${\cal S}_d(M)$ is connected, hence, by definition, ${\cal S}(M) = {\cal S}_d(M)$.
\qed

\vspace{2 mm}

\noindent
Department of Mathematics

\noindent
1 Shields Avenue

\noindent
University of
California, Davis

\noindent
Davis, CA 95616

\noindent
USA

\end{document}